\documentclass[a4paper,twoside,11pt]{article}
\usepackage{bbm}
\usepackage{amsfonts}
\usepackage{amssymb}\usepackage{amsmath}
 \numberwithin{equation}{section}
\begin{document}

\title{On the Cauchy problem of a weakly dissipative $\mu$HS equation}

\author{Jingjing Liu\footnote{e-mail:
jingjing830306@163.com}\ \ \ \ \ \ \
Zhaoyang Yin\footnote{e-mail: mcsyzy@mail.sysu.edu.cn}
\\ Department of Mathematics, Sun Yat-sen University,\\ 510275 Guangzhou, China}
\date{}
\maketitle

\begin{abstract}
In this paper, we study the Cauchy problem of a weakly dissipative $\mu$HS equation.
We first establish the local
well-posedness for the weakly dissipative $\mu$HS equation by Kato's semigroup theory. Then, we derive the precise
blow-up scenario for strong solutions to the equation. Moreover, we
present some blow-up results for strong solutions to the equation and give the blowup rate
of strong solutions to the equation when blowup occurs.
Finally, we give two global existence results to the equation.\\
\end{abstract}


\section{Introduction}
\par
Recently, the $\mu$-Hunter-Saxton (also called $\mu$-Camassa-Holm) equation
$$\mu(u)_{t}-u_{txx}=-2\mu(u)u_{x}+2u_{x}u_{xx}+uu_{xxx},$$
which is originally derived and studied in \cite{k-l-m} attracts a lot of attention. Here
$u(t,x)$ is a time-dependent function on the unit circle $\mathbb{S}=\mathbb{R}/\mathbb{Z}$ and
$\mu(u)=\int_{\mathbb{S}}udx$ denotes its mean. In \cite{k-l-m}, the authors show that if interactions
of rotators and an external magnetic field is allowed, then the $\mu$HS equation can
be viewed as a natural generalization of the rotator equation. Moreover, the $\mu$HS equation describes the geodesic flow
on $\mathcal{D}^{s}(\mathbb{S})$ with the right-invariant metric given at the
identity by the inner product \cite{k-l-m}
$$(u,v)=\mu(u)\mu(v)+\int_{\mathbb{S}}u_{x}v_{x}dx.$$ In \cite{k-l-m, l-m-t}, the authors showed that the $\mu$HS equation admits both periodic one-peakon solution and the multi-peakons. Moreover, in \cite{fu, gui}, the authors also discussed the $\mu$HS equation.

One of the closest relatives of the $\mu$HS equation is the Camassa-Holm equation
$$u_{t}-u_{txx}+3uu_{x}=2u_{x}u_{xx}+uu_{xxx},$$
which was introduced firstly by Fokas and Fuchssteiner in \cite{F-F} as
an abstract equation with a bihamiltonian structure. Meanwhile, it was derived by Camassa and Holm in \cite{R-D}
as a shallow water approximation independently. The Camassa-Holm equation is a model for shallow
water waves \cite{R-D, C-L, J-H, d} and a re-expression of the geodesic flow both
on the diffeomorphism group of the circle \cite{C-K} and on the Bott-Virasoro group
\cite{c-k-k-t, k}. The Camassa-Holm equation has a
bi-Hamiltonian structure \cite{F-F} and is completely integrable
\cite{c1, C-M}. The possibility of the relevance of Camassa-Holm to the modelling of tsunamis was
raised in  \cite{c-j, l}. It is worth to point
out that a long-standing open problem in hydrodynamics  was the derivation of a model equation
that can capture breaking waves as well as peaked traveling waves, cf. the discussion in \cite{w}.
The quest for peaked traveling waves is motivated by the desire to find waves replicating
a feature that is characteristic for the waves of great height-waves of largest amplitude that
are exact traveling solutions of the governing equations for water waves, whether periodic or solitary,
cf. \cite{c, t}. Breaking waves are solutions that remain bounded but their slope becomes unbounded in finite time, cf.\cite{c-e}. Both these aspects are modeled by the Camassa-Holm equation. Recently, the Camassa-Holm equation has been studied
extensively, c.f.\cite{rb, bc, bc1, c3, c4, c2, ce, ce1, cm, cs, dm, dm1, lo, m,  xz, y}. The other closest relatives of the $\mu$HS equation is the Hunter-Saxton equation \cite{J-R}
$$u_{txx}+2u_{x}u_{xx}+uu_{xxx}=0,$$
which is an asymptotic equation for rotators in liquid crystals and modeling the propagation of weakly nonlinear orientation waves in a
massive nematic liquid crystal. The orientation of the
molecules is described by the field of unit vectors $(\cos u(t,x),
\sin u(t,x))$ \cite{Y}. The Hunter-Saxton equation also arises in a
different physical context as the high-frequency limit \cite{HHD,
J-Y} of the Camassa-Holm equation. Similar to the Camassa-Holm equation, the Hunter-Saxton equation also has a bi-Hamiltonian
structure \cite{J-H, P-P} and is completely integrable
\cite{B1, J-Y}. The initial value problem of the Hunter-Saxton
equation also has been studied extensively, c.f.\cite{b-c, J-R, l1, Y}

 In general, it is difficult to avoid energy dissipation mechanisms in a real world. So, it is reasonable to study the model with energy dissipation. In \cite{g} and \cite{os}, the authors discussed the energy dissipative KdV equation from different aspects. Weakly dissipative CH equation and weakly dissipative DP equation have been studied in \cite{wy1} and \cite{ewy, wy, y1} respectively. Recently, some results for a weakly dissipative $\mu$DP equation are proved in \cite{k2}. It is worthy to note that the local well-posedness result in \cite{k2} is obtained by using a method based on a geometric argument.

 In this paper, we will discuss the Cauchy problem of the following weakly dissipative $\mu$HS equation:
\begin{equation}
\left\{\begin{array}{ll}
y_{t}+uy_{x}+2u_{x}y+\lambda y=0,&t > 0,\,x\in \mathbb{R},\\
y=\mu(u)-u_{xx},&t > 0,\,x\in \mathbb{R},\\
u(0,x) = u_{0}(x),& x\in \mathbb{R}, \\
u(t,x+1)=u(t,x), & t \geq 0, x\in \mathbb{R},\\\end{array}\right. \\
\end{equation}
or in the equivalent form:
\begin{equation}
\left\{\begin{array}{ll}
\mu(u)_{t}-u_{txx}+2\mu(u)u_{x}-2u_{x}u_{xx}-uu_{xxx}+\lambda(\mu(u)-u_{xx})=0,&t > 0,\,x\in \mathbb{R},\\
u(0,x) = u_{0}(x),& x\in \mathbb{R}, \\
u(t,x+1)=u(t,x), & t \geq 0, x\in \mathbb{R}.\\\end{array}\right. \\
\end{equation}
Here the constant $\lambda$ is assumed to be
positive and the term $\lambda y=\lambda(\mu(u)-u_{xx})$ models energy dissipation.

The paper is organized as follows. In Section 2, we establish the
local well-posedness of the initial value problem associated with
the equation (1.1). In Section 3, we derive the precise blow-up
scenario. In Section 4, we present some explosion criteria of strong
solutions to the equation (1.1) with general initial data and give the blowup rate
of strong solutions to the equation when blowup occurs. In Section
5, we give two new global existence results of strong solutions to the
equation (1.1).

 \textbf{Notation}  Given a Banach space $Z$, we denote its norm by
 $\|\cdot\|_{Z}$. Since all space of functions are over
 $\mathbb{S}=\mathbb{R}/\mathbb{Z}$, for simplicity, we drop $\mathbb{S}$ in our notations
 if there is no ambiguity. We let $[A,B]$ denote
 the commutator of linear operator $A$ and $B$. For convenience, we
 let $(\cdot|\cdot)_{s\times r}$ and $(\cdot|\cdot)_{s}$ denote the
 inner products of $H^{s}\times H^{r}$, $s,r\in \mathbb{R}_{+}$ and
 $H^{s}$, $s\in \mathbb{R}_{+}$, respectively.

\section{Local well-posedness}
\newtheorem {remark2}{Remark}[section]
\newtheorem{theorem2}{Theorem}[section]
\newtheorem{lemma2}{Lemma}[section]
In this section, we will establish the local well-posedness for the
Cauchy problem of the equation (1.1) in $H^{s}$, $s>\frac{3}{2}$, by applying Kato's theory.

\par
For convenience, we state here Kato's theory in the form suitable
for our purpose. Consider the abstract quasi-linear equation:
\begin{equation}
\frac{dv}{dt}+A(v)v=f(v),\ \ t>0,\ \ v(0)=v_{0}.
\end{equation}

\par
Let $X$ and $Y$ be Hilbert spaces such that $Y$ is continuously and
densely embedded in $X$ and let $Q : Y \rightarrow X$ be a
topological isomorphism. Let $L(Y,X)$ denote the space of all
bounded linear operator from $Y$ to $X$
($L(X)$, if $X=Y$.). Assume that: \\
(i) $A(y)\in L(Y,X)$ for $y\in Y$ with
$$ \|(A(y)-A(z))w\|_{X}\leq\mu_{1}\|
y-z\|_{X}\|w\|_{Y},\ \ \ \ y,z,w\in Y,$$ and $A(y)\in G(X,1,\beta)$,
(i.e. $A(y)$ is quasi-m-accretive),
uniformly on bounded sets in $Y$.\\
(ii) $QA(y)Q^{-1}=A(y)+B(y)$, where $B(y)\in L(X)$ is bounded,
uniformly on bounded sets in $Y$. Moreover,
$$ \|(B(y)-B(z))w\|_{X}\leq\mu_{2}\|
y-z\|_{Y}\|w\|_{X},\ \ \ \ y,z\in Y , \ w\in X.$$ (iii) $f:
Y\rightarrow Y$ and extends also to a map from $X$ to $X$. $f$ is
bounded on bounded sets in $Y$, and
$$ \|f(y)-f(z)\|_{Y}\leq \mu_{3}\|
y-z\|_{Y}, \ \ \ y,z\in Y,$$
$$
\|f(y)-f(z)\|_{X}\leq\mu_{4}\|y-z\|_{X}, \ \ \ y,z\in Y.$$ Here
$\mu_{1}, \mu_{2}, \mu_{3}$ and $ \mu_{4}$ depend only on
max$\{\|y\|_{Y}, \|z\|_{Y} \}$.
\begin{theorem2} \cite{Kato 1} Assume that (i), (ii) and (iii)
hold. Given  $ v_{0}\in Y$,  there exist a maximal $T>0$
 depending only on $\parallel v_{0}\parallel_{Y}$ and a
unique solution $v$ to Eq.(2.1) such that
$$
v=v(\cdot,v_{0})\in C([0,T);Y)\cap C^{1}([0,T);X).
$$
Moreover, the map $v_{0}\rightarrow v(\cdot,v_{0})$ is continuous
from Y to
$$
 C([0,T);Y)\cap C^{1}([0,T);X).
$$
\end{theorem2}

On one hand, with $y=\mu(u)-u_{xx},$ the first equation in (1.2) takes the form of a quasi-linear evolution
equation of hyperbolic type:
\begin{equation}
u_{t}+uu_{x}=-\partial_{x}A^{-1}(2\mu(u)u+\frac{1}{2}u_{x}^{2})-\lambda u,
\end{equation}
where $A=\mu-\partial_{x}^{2}$ is an isomorphism between
$H^{s}$ and $H^{s-2}$ with the inverse
$v=A^{-1}w$ given explicitly by \cite{escher, k-l-m}
\begin{align}
v(x)&=(\frac{x^{2}}{2}-\frac{x}{2}+\frac{13}{12})\mu(w)+(x-\frac{1}{2})\int_{0}^{1}\int_{0}^{y}w(s)dsdy\\
\nonumber&-\int_{0}^{x}\int_{0}^{y}w(s)dsdy+\int_{0}^{1}\int_{0}^{y}\int_{0}^{s}w(r)drdsdy.
\end{align}
Since $A^{-1}$ and
$\partial_{x}$ commute, the following identities hold
\begin{equation}
A^{-1}\partial_{x}w(x)=(x-\frac{1}{2})\int_{0}^{1}w(x)dx-\int_{0}^{x}w(y)dy+\int_{0}^{1}\int_{0}^{x}w(y)dydx,
\end{equation}
and
\begin{equation}
A^{-1}\partial_{x}^{2}w(x)=-w(x)+\int_{0}^{1}w(x)dx.
\end{equation}
On the other hand, integrating both sides of the first equation in (1.2) with respect
to $x$ on $\mathbb{S}$, we obtain
$$\frac{d}{dt}\mu(u)=-\lambda \mu(u),$$ it follows that
\begin{equation}
\mu(u)=\mu(u_{0})e^{-\lambda t}:=\mu_{0}e^{-\lambda t},
\end{equation}
where
$$\mu_{0}:=\mu(u_{0})=\int_{\mathbb{S}}u_{0}(x)dx.$$
Combining  (2.2) and (2.6), the equation (1.2) can be rewrite as
\begin{equation}
\left\{\begin{array}{ll}
u_{t}+uu_{x}=-\partial_{x}A^{-1}(2\mu_{0}e^{-\lambda t}u+\frac{1}{2}u_{x}^{2})-\lambda u,&t > 0,\,x\in \mathbb{R},\\
u(0,x) = u_{0}(x),& x\in \mathbb{R}, \\
u(t,x+1)=u(t,x), & t \geq 0, x\in \mathbb{R}.\\\end{array}\right. \\
\end{equation}

The remainder of this section is devoted to the local well-posedness result. Firstly, we will give a useful lemma.

\begin{lemma2} \cite{Kato 1} Let r,t be real
numbers such that $-r<t\leq r$. Then
$$
\left\|fg \right\|_{H^{t}}\leq
c\left\|f\right\|_{H^{r}}\left\|g\right\|_{H^{t}},  \ \ \ \  if \
r>\frac{1}{2},
$$
$$
\left\|fg \right\|_{H^{t+r-\frac{1}{2}}}\leq
c\left\|f\right\|_{H^{r}}\left\|g\right\|_{H^{t}},  \ if \
r<\frac{1}{2},
$$
where c is a positive  constant depending on r, t.
\end{lemma2}

\begin{theorem2}
Given $u_{0}\in H^{s}$, $s>\frac{3}{2},$
then there exists a maximal $T=T(u_{0})>0$, and a unique solution
$u$ to (2.7) (or (1.1)) such that
$$u=u(\cdot,u_{0})\in C([0,T); H^{s})\cap C^{1}([0,T);H^{s-1}).$$
Moreover, the solution depends continuously on the initial data, i.e., the
mapping $$u_{0}\rightarrow u(\cdot,u_{0}): H^{s}\rightarrow C([0,T); H^{s})\cap C^{1}([0,T);H^{s-1})$$
is continuous.
\end{theorem2}

\textbf{Proof}
For $u \in H^{s}$, $s>\frac{3}{2}$, we define the operator $A(u)= u\partial_{x}$. Similar to Lemma 2.6 in \cite{y2}, we have that
       $A(u)$ belongs to $G(H^{s-1},1,\beta)$, that
       is, $-A(u)$ generates a $C_{0}$-semigroup $T(t)$ on $H^{s-1}$ and $\|T(t)\|_{L(H^{s-1})}\leq
       e^{t\beta}$ for all $t\geq 0$. Analogous to Lemma 2.7 in
       \cite{y2}, we get that $A(u)\in L(H^{s},H^{s-1})$ and
$$
\left\|(A(z)-A(y))w\right\|_{H^{s-1}}\leq
\mu_{1}\left\|z-y\right\|_{H^{s-1}}\left\|w\right\|_{H^{s}}, $$ for all $z, y, w\in H^{s}$.

Let $Q=\Lambda=(1-\partial_{x}^{2})^{\frac{1}{2}}.$ Define $B(z)=QA(z)Q^{-1}-A(z)$ for $z\in H^{s},
s>\frac{3}{2}.$ Similar to Lemma 2.8 in \cite{y2}, we deduce that
$B(z)\in L(H^{s-1})$ and
$$
\left\|(B(z)-B(y))w\right\|_{H^{s-1}}\leq
\mu_{2}\left\|z-y\right\|_{H^{s}}\left\|w\right\|_{H^{s-1}}, $$ for all $z, y\in
H^{s}$ and $w\in H^{s-1}$. Where $\mu_{1},\ \mu_{2}$ are positive
constants.

Set $$f(u)=-\partial_{x}(\mu-\partial_{x}^{2})^{-1}(2\mu_{0}e^{-\lambda t}u+\frac{1}{2}u_{x}^{2})-\lambda u.$$
Let $y,z \in H^{s},$ $s>\frac{3}{2}.$ Since $H^{s-1}$ is an Banach algebra, it follows that
\begin{align*}
\|f(y)-f(z)\|_{H^{s}}\leq \ &\|-\partial_{x}(\mu-\partial_{x}^{2})^{-1}(2\mu_{0}e^{-\lambda t}(y-z)+\frac{1}{2}(y_{x}^{2}-z_{x}^{2}))\|_{H^{s}}+\lambda\|y-z\|_{H^{s}}\\
\leq \ &\|2\mu_{0}e^{-\lambda t}(y-z)+\frac{1}{2}(y_{x}+z_{x})(y_{x}-z_{x})\|_{H^{s-1}}+\lambda\|y-z\|_{H^{s}}\\
\leq \ &2|\mu_{0}|\|y-z\|_{H^{s-1}}+\frac{1}{2}\|y_{x}+z_{x}\|_{H^{s-1}}\|y_{x}-z_{x}\|_{H^{s-1}}+\lambda\|y-z\|_{H^{s}}\\
\leq \ &(2|\mu_{0}|+\|y\|_{H^{s}}+\|z\|_{H^{s}}+\lambda)\|y-z\|_{H^{s}}.
\end{align*}
Furthermore, taking $z=0$ in the above inequality, we obtain that $f$ is bounded on bounded set in $H^{s}.$ Moreover,
\begin{align*}
\|f(y)-f(z)\|_{H^{s-1}}\leq \ &\|-\partial_{x}(\mu-\partial_{x}^{2})^{-1}(2\mu_{0}e^{-\lambda t}(y-z)+\frac{1}{2}(y_{x}^{2}-z_{x}^{2}))\|_{H^{s-1}}+\lambda\|y-z\|_{H^{s-1}}\\
\leq \ &\|2\mu_{0}e^{-\lambda t}(y-z)+\frac{1}{2}(y_{x}+z_{x})(y_{x}-z_{x})\|_{H^{s-2}}+\lambda\|y-z\|_{H^{s-1}}\\
\leq \ &2|\mu_{0}|\|y-z\|_{H^{s-1}}+\frac{c}{2}\|y_{x}+z_{x}\|_{H^{s-1}}\|y_{x}-z_{x}\|_{H^{s-2}}+\lambda\|y-z\|_{H^{s-1}}\\
\leq \ &(2|\mu_{0}|+c(\|y\|_{H^{s}}+\|z\|_{H^{s}})+\lambda)\|y-z\|_{H^{s-1}},
\end{align*}
here we applied Lemma 2.1 with $r=s-1$, $t=s-2.$
Set $Y=H^{s},\ X=H^{s-1}.$ It
is obvious that $Q$ is an isomorphism of $Y$ onto $X$. Applying Theorem 2.1,
we obtain the local well-posedness
of Eq.(1.1) in $H^{s}$, $s>\frac{3}{2}$, and $u\in
C([0,T); H^{s})\cap C^{1}([0,T);H^{s-1})$. This completes the proof of Theorem 2.2. $\Box$

\begin{remark2}
Similar to the proof of Theorem 2.3 in \cite{y2}, we have that the maximal time of existence $T>0$ in Theorem 2.2 is independent
of the Sobolev index $s>\frac{3}{2}.$
\end{remark2}

\section{The precise blow-up scenario}
\newtheorem {remark3}{Remark}[section]
\newtheorem{theorem3}{Theorem}[section]
\newtheorem{lemma3}{Lemma}[section]
\newtheorem{definition3}{Definition}[section]
\newtheorem{claim3}{Claim}[section]

In this section, we present the precise blow-up scenario for strong
solutions to the equation (1.1). We first recall the following lemmas.

\begin{lemma3}
\cite{Kato 3} If $r>0$, then $H^{r}\cap L^{\infty}$ is an algebra.
Moreover
$$\parallel fg\parallel_{H^r}\leq c(\parallel  f\parallel_{L^{\infty}}\parallel g\parallel_{H^r}+\parallel
  f\parallel_{H^r}\parallel g\parallel_{L^{\infty}}),
$$
where c is a constant depending only on r.
\end{lemma3}

\begin{lemma3}
\cite{Kato 3} If $r>0$, then
$$\parallel[\Lambda^{r},f]g\parallel_{L^{2}}\leq c(\parallel \partial_{x}f\parallel_{L^{\infty}}
\parallel\Lambda^{r-1}g\parallel_{L^2}+\parallel
\Lambda^r f\parallel_{L^2}\parallel g\parallel_{L^{\infty}}),
$$
where c is a constant depending only on r.
\end{lemma3}

\begin{lemma3} \cite{constantin, fu} If $f\in H^{1}(\mathbb{S})$ is such that
$\int_{\mathbb{S}}f(x)dx=0,$ then we have
$$\max\limits_{x\in\mathbb{S}}f^{2}(x)\leq
\frac{1}{12}\int_{\mathbb{S}}f_{x}^{2}(x)dx.$$
\end{lemma3}

Next we prove the following useful result on global existence of
solutions to (1.1).
\begin{theorem3} Let $u_{0}\in H^{s}$, $s>\frac{3}{2}$, be given and assume that T is
the maximal existence time of the corresponding solution $u$
to (2.7) with the initial data $u_{0}$. If there exists $M>0$ such
that
$$\| u_{x}(t,\cdot)\|_{L^{\infty}}\leq M,\ \ t\in[0,T),$$
then the $H^s$-norm of $u(t,\cdot)$ does not blow up on [0,T).
\end{theorem3}
\textbf{Proof}
Let $u$ be the solution to (2.7) with the initial data $u_{0}\in H^s,\ s>\frac{3}{2}$,
and let T be the maximal existence time of the corresponding
solution $u$, which is guaranteed by Theorem 2.2. Throughout this
proof, $c>0$ stands for a generic constant depending only on $s$.

Applying the operator $\Lambda^{s}$ to the first equation in (2.7),
multiplying by $\Lambda^{s} u$, and integrating over $\mathbb{S}$,
we obtain
\begin{equation}
\frac{d}{dt}\|u\|^{2}_{H^{s}}=-2(uu_{x},
u)_{s}-2(u, \partial_{x}(\mu-\partial_{x}^{2})^{-1}(2\mu_{0}e^{-\lambda t}u+\frac{1}{2}u_{x}^{2}))_{s}-2\lambda(u,u)_{s}.
\end{equation}
Let us estimate the first term of the right hand side of (3.1).
\begin{align}
\nonumber|(uu_{x},u)_{s}|&= |(\Lambda^s(u\partial_{x}u),\Lambda^s u)_{0}|\\
\nonumber&=|([\Lambda^s,u]\partial_{x}u,\Lambda^s
u)_{0}+(u\Lambda^s\partial_{x}u,\Lambda^s u)_{0}|\\
\nonumber&\leq \|[\Lambda^s,u]\partial_{x}u\|_{L^2}\|\Lambda^s
u\|_{L^2}+\frac{1}{2}|(u_{x}\Lambda^s u,\Lambda^s u)_{0}|\\
\nonumber& \leq (c\| u_{x}\|_{L^{\infty}}+\frac{1}{2}\|
u_{x}\|_{L^{\infty}})\|u\|^2_{H^s}\\
&\leq c\|u_{x}\|_{L^{\infty}}\|u\|^2_{H^s},
\end{align}
where we used Lemma 3.2 with $r=s$. Furthermore, we estimate the second term of the right
hand side of (3.1) in the following way:
\begin{align}
\nonumber&|(u, \partial_{x}(\mu-\partial_{x}^{2})^{-1}(2\mu_{0}e^{-\lambda t}u+\frac{1}{2}u_{x}^{2}))_{s}|\\
\nonumber\leq \ &\|\partial_{x}(\mu-\partial_{x}^{2})^{-1}(2\mu_{0}e^{-\lambda t}u+\frac{1}{2}u_{x}^{2})\|_{H^{s}}\|u\|_{H^{s}}\\
\nonumber\leq \ &\|2\mu_{0}e^{-\lambda t}u+\frac{1}{2}u_{x}^{2}\|_{H^{s-1}}\|u\|_{H^{s}}\\
\nonumber\leq \ &c(|\mu_{0}|\|u\|_{H^{s}}+\|u_{x}\|_{L^{\infty}}\|u_{x}\|_{H^{s-1}})\|u\|_{H^{s}}\\
\leq \ &c(|\mu_{0}|+\|u_{x}\|_{L^{\infty}})\|u\|_{H^{s}}^{2},
\end{align}
where we used Lemma 3.1 with $r=s-1.$ Combining (3.2) and (3.3) with
(3.1), we get
$$\frac{d}{dt}\|u\|^{2}_{H^{s}}\leq c(|\mu_{0}|+\|u_{x}\|_{L^{\infty}}+2\lambda)\|u\|_{H^{s}}^{2}.$$
An application of Gronwall's inequality and the assumption of the
theorem yield
$$\|u\|_{H^{s}}^{2}\leq e^{c(|\mu_{0}|+M+2\lambda)t}\|u_{0}\|_{H^{s}}^{2}.$$
This completes the proof of the theorem. $\Box$ \\

The following result describes the precise blow-up scenario for sufficiently regular solutions to the equation (1.1).
\begin{theorem3}
Let $u_{0}\in H^{s}$, $s>\frac{3}{2} $ be given and let T be the
maximal existence time of the corresponding solution $u$ to (2.7) with the initial data $u_{0}$. Then the corresponding
solution blows up in finite time if and only if
$$\liminf\limits_{t\rightarrow T}\{\inf\limits_{x\in\mathbb{S}}u_{x}(t,x)\}=-\infty.$$
\end{theorem3}
\textbf{Proof} Applying a simple density argument, Remark 2.1 implies that we only need to consider the case $s=3.$ Multiplying the first equation
in (1.1) by $y$ and integrating over $\mathbb{S}$ with respect to $x$ yield
\begin{align*}
\frac{d}{dt}\int_{\mathbb{S}}y^{2}dx= \ &2\int_{\mathbb{S}}y(-uy_{x}-2u_{x}y-\lambda y)dx\\
= \ &-2\int_{\mathbb{S}}uyy_{x}dx-4\int_{\mathbb{S}}u_{x}y^{2}dx-2\lambda\int_{\mathbb{S}}y^{2}dx\\
= \ &-3\int_{\mathbb{S}}u_{x}y^{2}dx-2\lambda\int_{\mathbb{S}}y^{2}dx.
\end{align*}
So, if there is a constant $M>0$ such that $$u_{x}(t,x)\geq -M, \ \ \ \ \ \forall \ (t,x)\in[0,T)\times \mathbb{S},$$
then $$\frac{d}{dt}\int_{\mathbb{S}}y^{2}dx \leq (3M-2\lambda)\int_{\mathbb{S}}y^{2}dx.$$
Gronwall's inequality implies that
$$
\int_{\mathbb{S}}y^{2}dx\leq e^{(3M-2\lambda)}\int_{\mathbb{S}}y^{2}(0,x)dx.
$$
Note that
$$\int_{\mathbb{S}}y^{2}dx=\mu(u)^{2}+\int_{\mathbb{S}}u_{xx}^{2}dx\geq \|u_{xx}\|_{L^{2}}^{2}.$$
Since $u_{x}\in H^{2}\subset H^{1}$ and $\int_{\mathbb{S}}u_{x}dx=0,$ Lemma 3.3 implies that
$$\|u_{x}\|_{L^{\infty}}\leq \frac{1}{2\sqrt{3}}\|u_{xx}\|_{L^{2}}\leq e^{\frac{3M-2\lambda}{2}}\|y(0,x)\|_{L^{2}}.$$
Theorem 3.1 ensures that the solution $u$ does not blow up in finite time.

On the other hand, by Sobolev's imbedding theorem it is clear that if
$$\liminf\limits_{t\rightarrow T}\{\inf\limits_{x\in\mathbb{S}}u_{x}(t,x)\}=-\infty,$$ then $T<\infty$.
This completes the proof of the theorem. $\Box$

\section{Blow-up and blow-up rate}
\newtheorem{theorem4}{Theorem}[section]
\newtheorem{lemma4}{Lemma}[section]
\newtheorem {remark4}{Remark}[section]
\newtheorem{corollary4}{Corollary}[section]

In this section, we discuss the blow-up phenomena of the equation
(1.1) and prove that there exist strong solutions to (1.1) which do
not exist globally in time.

Firstly, for $u_{0}\in H^{s}, \  s>\frac{3}{2},$ we will give some useful estimates for the corresponding solution $u.$
By the first equation of (1.2) and (2.6), a direct computation implies that
\begin{align*}
\frac{d}{dt}\int_{\mathbb{S}}u_{x}^{2}dx= \ &2\int_{\mathbb{S}}u(-u_{txx})dx\\
= \ &2\int_{\mathbb{S}}u(-\mu(u)_{t}-2\mu(u)u_{x}+2u_{x}u_{xx}+uu_{xxx}-\lambda\mu(u)+\lambda u_{xx})dx\\
= \ &-2\mu(u)_{t}\mu(u)-2\lambda(\mu(u))^{2}-2\lambda\int_{\mathbb{S}}u_{x}^{2}dx\\
= \ &-2\lambda\int_{\mathbb{S}}u_{x}^{2}dx.
\end{align*}
It follows that
\begin{equation}
\int_{\mathbb{S}}u_{x}^{2}dx=\int_{\mathbb{S}}u_{0,x}^{2}dx\cdot e^{-2\lambda t}:=\mu_{1}^{2}e^{-2\lambda t},
\end{equation}
where $\mu_{1}=\left(\int_{\mathbb{S}}u_{0,x}^{2}dx\right)^{\frac{1}{2}}.$
Note that $\int_{\mathbb{S}}(u(t,x)-\mu(u))dx=\mu(u)-\mu(u)=0.$
By Lemma 3.3, we find that
$$\max\limits_{x\in\mathbb{S}}[u(t,x)-\mu(u)]^{2}\leq
 \frac{1}{12}\int_{\mathbb{S}}u_{x}^{2}(t,x)dx\leq
 \frac{1}{12}\mu_{1}^{2}.$$ So we have
 \begin{equation}
 \|u(t,\cdot)\|_{L^{\infty}}\leq
 |\mu_{0}|+\frac{\sqrt{3}}{6}\mu_{1}.
 \end{equation}

\begin{lemma4}
\cite{c-e}
Let $t_{0}>0$ and $v\in C^{1}([0,t_{0}); H^{2}(\mathbb{R}))$. Then
for every $t\in[0,t_{0})$ there exists at least one point $\xi(t)\in
\mathbb{R}$ with
$$ m(t):=\inf_{x\in \mathbb{R}}\{v_{x}(t,x)\}=v_{x}(t,\xi(t)),$$ and
the function $m$ is almost everywhere differentiable on $(0,t_{0})$
with $$ \frac{d}{dt}m(t)=v_{tx}(t,\xi(t)) \ \ \ \ a.e.\ on \
(0,t_{0}).$$
\end{lemma4}

\begin{theorem4}Let $u_{0}\in H^s, s>\frac{3}{2},$ $u_{0}\not\equiv c$ for $\forall \ c\in \mathbb{R}$ and T be the maximal time of the
solution $u$ to (1.1) with the initial data $u_0$. If $u_{0}$ satisfies the following condition
$$\int_{\mathbb{S}}u_{0,x}^{3}dx<-3\lambda \mu_{1}^{2}-\mu_{1}\sqrt{9\lambda^{2}\mu_{1}^{2}+2K},$$ where
$K=6|\mu_{0}|\mu_{1}^{2}(|\mu_{0}|+\frac{\sqrt{3}}{6}\mu_{1}),$
then the corresponding solution to (1.1) blows up in finite time.
\end{theorem4}
\textbf{Proof} As mentioned earlier, here we only need
to show that the above theorem holds for $s=3$. Differentiating the
first equation of Eq.(2.7) with respect to $x$, we have
\begin{equation}
u_{tx}=-\frac{1}{2}u_{x}^{2}-uu_{xx}+2\mu_{0}e^{-\lambda t}u-\lambda u_{x}-2\mu_{0}^{2}e^{-2\lambda t}-\frac{1}{2}\mu_{1}^{2}e^{-2\lambda t}
\end{equation}
Then, it follows that
\begin{align*}
&\frac{d}{dt}\int_{\mathbb{S}}u_{x}^{3}dx\\
=&\int_{\mathbb{S}}3u_{x}^{2}u_{xt}dx\\
=&3\int_{\mathbb{S}}u_{x}^{2}(-\frac{1}{2}u_{x}^{2}-uu_{xx}+2\mu_{0}e^{-\lambda t}u-\lambda u_{x}-2\mu_{0}^{2}e^{-2\lambda t}-\frac{1}{2}\mu_{1}^{2}e^{-2\lambda t})dx\\
\leq & -\frac{3}{2}\int_{\mathbb{S}}u_{x}^{4}dx-3\int_{\mathbb{S}}uu_{x}^{2}u_{xx}dx+6\mu_{0}e^{-\lambda t}\int_{\mathbb{S}}uu_{x}^{2}dx-3\lambda\int_{\mathbb{S}}u_{x}^{3}dx\\
=&-\frac{1}{2}\int_{\mathbb{S}}u_{x}^{4}dx+6\mu_{0}e^{-\lambda t}\int_{\mathbb{S}}uu_{x}^{2}dx-3\lambda\int_{\mathbb{S}}u_{x}^{3}dx\\
\leq &-\frac{1}{2}\int_{\mathbb{S}}u_{x}^{4}dx-3\lambda\int_{\mathbb{S}}u_{x}^{3}dx+6|\mu_{0}|\mu_{1}^{2}(|\mu_{0}|+\frac{\sqrt{3}}{6}\mu_{1})\\
:= & -\frac{1}{2}\int_{\mathbb{S}}u_{x}^{4}dx-3\lambda\int_{\mathbb{S}}u_{x}^{3}dx+K.
\end{align*}
Using the following inequality
$$
\left|\int_{\mathbb{S}}u_{x}^{3}dx\right| \leq
\left(\int_{\mathbb{S}}u_{x}^{4}dx\right)^{\frac{1}{2}}\left(\int_{\mathbb{S}}u_{x}^{2}dx\right)^{\frac{1}{2}}
\leq\left(\int_{\mathbb{S}}u_{x}^{4}dx\right)^{\frac{1}{2}}\mu_{1},$$ and letting
$$m(t)=\int_{\mathbb{S}}u_{x}^{3}dx,$$ we have
\begin{align*}
\frac{d}{dt}m(t)\leq & -\frac{1}{2\mu_{1}^{2}}m^{2}(t)-3\lambda m(t)+K\\
= &-\frac{1}{2\mu_{1}^{2}}\left(m(t)+3\lambda \mu_{1}^{2}+\mu_{1}\sqrt{9\lambda^{2}\mu_{1}^{2}+2K}\right)\left(m(t)+3\lambda \mu_{1}^{2}-\mu_{1}\sqrt{9\lambda^{2}\mu_{1}^{2}+2K}\right).
\end{align*}
Note
that if $m(0)<-3\lambda \mu_{1}^{2}-\mu_{1}\sqrt{9\lambda^{2}\mu_{1}^{2}+2K}$ then
$m(t)<-3\lambda \mu_{1}^{2}-\mu_{1}\sqrt{9\lambda^{2}\mu_{1}^{2}+2K}$ for all
$t\in[0,T)$. From the above inequality we obtain
$$\frac{m(0)+A+B}{m(0)+A-B}e^{\frac{B}{\mu_{1}^{2}}t}-1\leq \frac{2B}{m(t)+A-B}\leq 0$$
with $A=3\lambda \mu_{1}^{2},$ $B=\mu_{1}\sqrt{9\lambda^{2}\mu_{1}^{2}+2K}.$ Since $0<\frac{m(0)+A+B}{m(0)+A-B}<1$,
then there exists
$$0<T\leq \frac{\mu_{1}}{\sqrt{9\lambda^{2}\mu_{1}^{2}+2K}}\ln\frac{m(0)+3\lambda \mu_{1}^{2}-\mu_{1}\sqrt{9\lambda^{2}\mu_{1}^{2}+2K}}{m(0)+3\lambda \mu_{1}^{2}+\mu_{1}\sqrt{9\lambda^{2}\mu_{1}^{2}+2K}},$$
such that $\lim_{t\rightarrow T} m(t)=-\infty.$ On the other hand,
$$\int_{\mathbb{S}}u_{x}^{3}dx\geq\inf\limits_{x\in\mathbb{S}}u_{x}(t,x)\int_{\mathbb{S}}u_{x}^{2}dx=\inf\limits_{x\in\mathbb{S}}u_{x}(t,x)\cdot
\mu_{1}^{2}e^{-2\lambda t}.$$
Applying Theorem 3.2, the solution $u$ blows up in finite time. $\Box$

\begin{theorem4} Let $u_{0}\in H^s, s>\frac{3}{2},$ and T be the maximal time of the
solution $u$ to (1.1) with the initial data $u_0$. If
$$\inf\limits_{x\in \mathbb{S}}u_{0}^{\prime}(x)<-\lambda-\sqrt{\lambda^{2}+2K},$$ with
$K=2|\mu_{0}|(|\mu_{0}|+\frac{\sqrt{3}}{6}\mu_{1}),$
then the corresponding solution to (1.1) blows up in finite time.
\end{theorem4}
\textbf{Proof} As mentioned earlier, here we only need
to show that the above theorem holds for $s=3$. Define now $$ m(t):=\min\limits_{x\in \mathbb{S}}[u_{x}(t,x)], \ \ t\in
[0,T)$$ and let $\xi(t)\in \mathbb{S}$ be a point where this
minimum is attained by using Lemma 4.1. It follows that
$$m(t)=u_{x}(t,\xi(t)).$$
Clearly $u_{xx}(t,\xi(t))=0$ since $u(t,\cdot)\in
H^{3}(\mathbb{S})\subset C^{2}(\mathbb{S}).$ Evaluating (4.3) at
$(t,\xi(t))$, we obtain
\begin{align*}
\frac{dm(t)}{dt}=&-\frac{1}{2}m^{2}(t)+2\mu_{0}e^{-\lambda t}u(t,\xi(t))-\lambda m(t)-2\mu_{0}^{2}e^{-2\lambda t}-\frac{1}{2}\mu_{1}^{2}e^{-2\lambda t}\\
\leq & -\frac{1}{2}m^{2}(t)-\lambda m(t)+2|\mu_{0}|(|\mu_{0}|+\frac{\sqrt{3}}{6}\mu_{1})\\
:= & -\frac{1}{2}m^{2}(t)-\lambda m(t)+K\\
= & -\frac{1}{2}(m(t)+\lambda+\sqrt{\lambda^{2}+2K})(m(t)+\lambda-\sqrt{\lambda^{2}+2K}).
\end{align*}
Note
that if $m(0)<-\lambda-\sqrt{\lambda^{2}+2K}$ then $m(t)<-\lambda-\sqrt{\lambda^{2}+2K}$ for all
$t\in[0,T)$. From the above inequality we obtain
$$\frac{m(0)+\lambda+\sqrt{\lambda^{2}+2K}}{m(0)+\lambda-\sqrt{\lambda^{2}+2K}}e^{\sqrt{\lambda^{2}+2K}t}-1\leq \frac{2\sqrt{\lambda^{2}+2K}}{m(t)+\lambda-\sqrt{\lambda^{2}+2K}}\leq 0$$
Since $0<\frac{m(0)+\lambda+\sqrt{\lambda^{2}+2K}}{m(0)+\lambda-\sqrt{\lambda^{2}+2K}}<1$,
then there exists
$$0<T\leq \frac{1}{\sqrt{\lambda^{2}+2K}}\ln\frac{m(0)+\lambda-\sqrt{\lambda^{2}+2K}}{m(0)+\lambda+\sqrt{\lambda^{2}+2K}},$$
such that $\lim_{t\rightarrow T} m(t)=-\infty.$ Theorem 3.2 implies the solution $u$ blows up in finite time. $\Box$

\begin{theorem4} Let $u_{0}\in H^s, s>\frac{3}{2},$ and T be the maximal time of the
solution $u$ to (1.1) with the initial data $u_0$. If $u_{0}(x)$ is odd satisfies $u_{0}^{\prime}(0)<-2\lambda,$
then the corresponding solution to (1.1) blows up in finite time.
\end{theorem4}
\textbf{Proof} As mentioned earlier, here we only need
to show that the above theorem holds for $s=3$. By $\mu(-u(t,-x))=-\mu(u(t,x)),$ we have (1.2) is invariant under the transformation
$(u,x)\rightarrow(-u,-x).$ Thus we deduce that if $u_{0}(x)$ is odd, then $u(t,x)$ is odd with respect to $x$ for any
$t\in[0,T).$ By continuity with respect to $x$ of $u$ and $u_{xx},$
we have$$
u(t,0)=u_{xx}(t,0)=0, \ \ \ \forall \ t\in[0,T).$$ Evaluating (4.3) at
$(t,0)$ and letting $h(t)=u_{x}(t,0),$ we obtain
\begin{align*}
\frac{dh(t)}{dt}=&-\frac{1}{2}h^{2}(t)-\lambda h(t)-2\mu_{0}^{2}e^{-2\lambda t}-\frac{1}{2}\mu_{1}^{2}e^{-2\lambda t}\\
\leq & -\frac{1}{2}h^{2}(t)-\lambda h(t)\\
= & -\frac{1}{2}h(t)(h(t)+2\lambda).
\end{align*}
Note
that if $h(0)<-2\lambda$ then $h(t)<-2\lambda$ for all
$t\in[0,T)$. From the above inequality we obtain
$$\left(1+\frac{2\lambda}{h(0)}\right)e^{\lambda t}-1\leq\frac{2\lambda}{h(t)}\leq 0.$$
Since $\frac{h(0)}{h(0)+2\lambda}> 1,$ then there exists
$$0<T\leq\frac{1}{\lambda}\ln\frac{h(0)}{h(0)+2\lambda},$$
such that $\lim_{t\rightarrow T} h(t)=-\infty.$ Theorem 3.2 implies the solution $u$ blows up in finite time. $\Box$\\

Consequently, we will discuss the blow-up rate for the wave-breaking
solutions to Eq.(1.1). The following result implies that the blow-up rate of
strong solutions to weakly dissipative $\mu$-HS equation is not affected by the weakly dissipative term even though the occurrence of blow-up of strong solutions to Eq. (1.1) is affected by
the dissipative parameter, see Theorem 4.1-4.3.

\begin{theorem4}
Let $u_{0}\in H^s, s>\frac{3}{2},$ and T be the maximal time of the
solution $u$ to (1.1) with the initial data $u_0$. If $T$ is finite, we obtain
$$ \lim_{t\rightarrow T}(T-t)\min_{x\in \mathbb{S}}u_{x}(t,x)=-2.$$
\end{theorem4}
\textbf{Proof} From the proof of Theorem 4.2, with $ m(t):=\min\limits_{x\in \mathbb{S}}[u_{x}(t,x)], \ t\in
[0,T),$ we have
\begin{align*}\left|m^{\prime}(t)+\frac{1}{2}m^{2}(t)+\lambda m(t)\right|=&\left|2\mu_{0}e^{-\lambda t}u(t,\xi(t))-2\mu_{0}^{2}e^{-2\lambda t}-\frac{1}{2}\mu_{1}^{2}e^{-2\lambda t}\right|\\
\leq & 2|\mu_{0}|(|\mu_{0}|+\frac{\sqrt{3}}{6}\mu_{1})+2\mu_{0}^{2}+\frac{1}{2}\mu_{1}^{2}:=K
\end{align*}
It follows that
\begin{equation}
-K\leq m^{\prime}(t)+\frac{1}{2}m^{2}(t)+\lambda m(t)\leq K \ \ \ a.e.\ \ on\ (0,T).
\end{equation}
Thus,
\begin{equation}
-K-\frac{1}{2}\lambda^{2}\leq m^{\prime}(t)+\frac{1}{2}(m(t)+\lambda)^{2}\leq K+\frac{1}{2}\lambda^{2}\ \ \ a.e.\ \ on\ (0,T).
\end{equation}

Let $\varepsilon\in(0,\frac{1}{2})$. Since $\liminf\limits_{t\rightarrow
T}m(t)= -\infty$ by Theorem 3.2, there is some $t_{0}\in (0,T)$ with
$m(t_{0})+\lambda<0$ and $(m(t_{0})+\lambda)^{2}>\frac{1}{\varepsilon}(K+\frac{1}{2}\lambda^{2})$. Let us first
prove that
\begin{equation}
(m(t)+\lambda)^{2}>\frac{1}{\varepsilon}(K+\frac{1}{2}\lambda^{2}), \ \ \ \ t\in [t_{0}, T).
\end{equation}

Since $m$ is locally Lipschitz (it belongs to
$W_{loc}^{1,\infty}(\mathbb{R})$ by Lemma 4.1) there is some
$\delta>0$ such that
$$(m(t)+\lambda)^{2}>\frac{1}{\varepsilon}(K+\frac{1}{2}\lambda^{2}), \ \ \ \ t\in [t_{0},
t_{0}+\delta).$$ Pick $\delta>0$ maximal with this property. If
$\delta < T-t_{0}$ we would have
$(m(t_{0}+\delta)+\lambda)^{2}=\frac{1}{\varepsilon}(K+\frac{1}{2}\lambda^{2})$ while$$
 m^{\prime}(t)\leq-\frac{1}{2}(m(t)+\lambda)^{2}+ K+\frac{1}{2}\lambda^{2}<(\varepsilon-\frac{1}{2})(m(t)+\lambda)^{2}< 0 \ \ \ \ \ \ a.e.\ on\ (t_{0}, t_{0}+\delta).$$ Being locally
Lipschitz, the function $m$ is absolutely continuous and therefore
we would obtain by integrating the previous relation on $[t_{0},
t_{0}+\delta]$ that $$ m(t_{0}+\delta)+\lambda\leq m(t_{0})+\lambda<0,$$ which on
its turn would yield$$(m(t_{0}+\delta)+\lambda)^{2}\geq
(m(t_{0})+\lambda)^{2}>\frac{1}{\varepsilon}(K+\frac{1}{2}\lambda^{2}).$$ The obtained contradiction
completes the proof of the relation (4.6).

A combination of (4.5) and (4.6) enables us to infer
\begin{equation}
\frac{1}{2}+\varepsilon\geq-\frac{m^{\prime}(t)}{(m(t)+\lambda)^{2}}\geq\frac{1}{2}-\varepsilon
\ \ \ a.e.\ on \ (0,T).
\end{equation}
Since $m(t)+\lambda$ is locally Lipschitz on $[0,T)$ and (4.6) holds, it is
easy to check that $\frac{1}{m(t)+\lambda}$ is locally Lipschitz on
$(t_{0},T).$ Differentiating the relation $(m(t)+\lambda)\cdot
\frac{1}{m(t)+\lambda}=1,\ t\in\ (t_{0}, T),$ we get
$$\left(\frac{1}{m(t)+\lambda}\right)^{\prime}=-\frac{m^{\prime}(t)}{(m(t)+\lambda)^{2}} \ a.e.\ on\
(t_{0}, T),$$ with $\frac{1}{m(t)+\lambda}$ absolutely continuous on
$(t_{0},T).$ For $t\in (t_{0},T)$. Integrating (4.7) on $(t,T)$ to
obtain $$
(\frac{1}{2}+\varepsilon)(T-t)\geq-\frac{1}{m(t)+\lambda}\geq(\frac{1}{2}-\varepsilon)(T-t),\
t\in (t_{0}, T),$$ that is,$$
\frac{1}{\frac{1}{2}+\varepsilon}\leq-(m(t)+\lambda)(T-t)\leq\frac{1}{\frac{1}{2}-\varepsilon},\
t\in(t_{0},T).$$ By the arbitrariness of
$\varepsilon\in(0,\frac{1}{2}),$ we have
$$\lim\limits_{t\rightarrow T}(m(t)+\lambda)(T-t)=-2,$$ so
the statement of Theorem 4.4
follows.  $\Box$

\section{Global Existence}
\newtheorem{theorem5}{Theorem}[section]
\newtheorem{lemma5}{Lemma}[section]
\newtheorem {remark5}{Remark}[section]
\newtheorem{corollary5}{Corollary}[section]

In this section, we will present some global existence results. Firstly,
we give a useful lemma.

Given $u_{0}\in H^{s}$ with $s>\frac{3}{2}$. Theorem 2.2
ensures the existence of a maximal  $T> 0$ and a solution $u$ to (2.7) such that
$$u=u(\cdot,u_{0})\in C([0,T); H^{s})\cap C^{1}([0,T);H^{s-1}).$$
Consider now the following initial value problem
\begin{equation}
\left\{\begin{array}{ll}q_{t}=u(t,q),\ \ \ \ t\in[0,T), \\
q(0,x)=x,\ \ \ \ x\in\mathbb{R}. \end{array}\right.
\end{equation}

\begin{lemma5}
Let $u_{0}\in H^{s}$ with $s>\frac{3}{2},$ $T> 0$ be the maximal existence time. Then Eq.(5.1) has a unique solution
$q\in C^1([0,T)\times \mathbb{R};\mathbb{R})$ and the map
$q(t,\cdot)$ is an increasing diffeomorphism of $\mathbb{R}$ with
$$
q_{x}(t,x)=exp\left(\int_{0}^{t}u_{x}(s,q(s,x))ds\right)>0, \ \
(t,x)\in [0,T)\times \mathbb{R}.$$
Moreover, with $y=\mu(u)-u_{xx},$ we have
$$
y(t,q(t,x))q_{x}^{2}=y_{0}(x)e^{-\lambda t}.
$$
\end{lemma5}
\textbf{Proof}
The proof of the first conclusion is similar to the proof of Lemma 4.1 in \cite{Y1}, so we omit it here. By the first equation in (1.1) and the equation (5.1), we have
\begin{align*}
&\frac{d}{dt}y(t,q(t,x))q_{x}^{2}\\
= \ &(y_{t}+y_{x}q_{t})q_{x}^{2}+y\cdot2q_{x}q_{xt}\\
= \ &(y_{t}+uy_{x})q_{x}^{2}+2yu_{x}q_{x}^{2}\\
= \ &(y_{t}+uy_{x}+2yu_{x})q_{x}^{2}=-\lambda y q_{x}^{2}.
\end{align*}
It follows that $$y(t,q(t,x))q_{x}^{2}=y_{0}(x)e^{-\lambda t}. \ \ \ \ \ \ \ \ \ \ \ \ \ \ \ \Box$$

\begin{theorem5}
If $y_{0}=\mu_{0}-u_{0,xx}\in H^{1}$ does not change sign, then the corresponding solution $u$ of the initial value $u_{0}$ exists globally in time.
\end{theorem5}
\textbf{Proof}
Note that given $t\in[0,T),$ there is a $\xi(t)\in \mathbb{S}$ such that $u_{x}(t, \xi(t))=0$ by the periodicity of $u$ to $x$-variable. If $y_{0}\geq 0,$ then Lemma 5.1 implies that $y\geq 0.$ For $x\in[\xi(t), \xi(t)+1],$ we have
\begin{align*}
-u_{x}(t,x)=-\int_{\xi(t)}^{x}\partial_{x}^{2}u(t,x)dx=&\int_{\xi(t)}^{x}(y-\mu(u))dx=\int_{\xi(t)}^{x}ydx-\mu(u)(x-\xi(t))\\
\leq \ &\int_{\mathbb{S}}ydx-\mu(u)(x-\xi(t))=\mu(u)(1-x+\xi(t))\leq |\mu_{0}|.
\end{align*}
It follows that $u_{x}(t,x)\geq -|\mu_{0}|.$
 On the other hand, if $y_{0}\leq 0,$ then Lemma 5.1 implies that $y\leq 0.$ Therefore, for $x\in[\xi(t), \xi(t)+1],$ we have
\begin{align*}
-u_{x}(t,x)=-\int_{\xi(t)}^{x}\partial_{x}^{2}u(t,x)dx=&\int_{\xi(t)}^{x}(y-\mu(u))dx=\int_{\xi(t)}^{x}ydx-\mu(u)(x-\xi(t))\\
\leq \ &-\mu(u)(x-\xi(t))\leq |\mu_{0}|.
\end{align*}
It follows that $u_{x}(t,x)\geq -|\mu_{0}|.$
This completes the proof by using Theorem 3.3. $\Box$

\begin{corollary5}
If the initial value $u_{0}\in H^{3}$ such that $$\|\partial_{x}^{3}u_{0}\|_{L^{2}}\leq 2\sqrt{3}|\mu_{0}|,$$ then the corresponding solution $u$ of $u_{0}$ exists globally in time.
\end{corollary5}
\textbf{Proof} Note that $\int_{\mathbb{S}}\partial_{x}^{2}u_{0}dx=0,$ Lemma 3.3 implies that $$\|\partial_{x}^{2}u_{0}\|_{L^{\infty}}\leq \frac{\sqrt{3}}{6}\|\partial_{x}^{3}u_{0}\|_{L^{2}}.$$
If $\mu_{0}\geq 0$, then
$$y_{0}=\mu_{0}-\partial_{x}^{2}u_{0}\geq \mu_{0}-\frac{\sqrt{3}}{6}\|\partial_{x}^{3}u_{0}\|_{L^{2}}\geq\mu_{0}-|\mu_{0}|=0.$$
If $\mu_{0}\leq 0$, then
$$y_{0}=\mu_{0}-\partial_{x}^{2}u_{0}\leq \mu_{0}+\|\partial_{x}^{2}u_{0}\|_{L^{\infty}}\leq\mu_{0}+\frac{\sqrt{3}}{6}\|\partial_{x}^{3}u_{0}\|_{L^{2}}\leq\mu_{0}+|\mu_{0}|=0$$
This completes the proof by using Theorem 5.1. $\Box$

\bigskip

\noindent\textbf{Acknowledgments} This work was partially supported by NNSFC (No. 10971235), RFDP (No. 200805580014), NCET (No. 08-0579) and the key project of Sun Yat-sen University. The authors thank the
referees for useful comments and suggestions.


\begin{thebibliography}{99}
\small

\bibitem{B1}R. Beals, D. Sattinger and J. Szmigielski, Inverse scattering solutions of the Hunter-Saxton
equations, {\it Appl. Anal.}, {\bf78} (2001), 255-269.

\bibitem{rb}G. Blanco, On the Cauchy problem for the Camassa-Holm equation, {\it Nonlinear
Anal.}, {\bf 46} (2001), 309-327.

\bibitem{b-c}A. Bressan and A. Constantin, Global solutions of the Hunter-Saxton equation, {\it SIAM J. Math. Anal.},
{\bf 37} (2005), 996-1026.

\bibitem{bc}A. Bressan and A. Constantin, Global conservative solutions of the Camassa-Holm equation,
{\it Arch. Rat. Mech. Anal.}, {\bf 183} (2007), 215-239.

\bibitem{bc1}A. Bressan and A. Constantin, Global dissipative solutions of the Camassa-Holm equation,
{\it Anal. Appl.}, {\bf 5} (2007), 1-27.

\bibitem{R-D}R. Camassa and D. Holm, An integrable shallow water equation with peaked solitons, {\it Phys.
Rev. Lett.}, {\bf 71} (1993), 1661-1664.

\bibitem{c3}A. Constantin, On the Cauchy problem for the periodic Camassa-Holm equation, {\it J. Differential
Equations}, {\bf 141} (1997), 218-235.

\bibitem{c1}A. Constantin, On the inverse spectral problem for the
Camassa-Holm equation, {\it J. Funct. Anal.}, {\bf 155} (1998), 352-363.

\bibitem{c4}A. Constantin, On the Blow-up of solutions of a periodic shallow water equation, {\it J. Nonlinear
Sci.}, {\bf 10} (2000), 391-399.

\bibitem{c2}A. Constantin, Existence of permanent and breaking
waves for a shallow water equation: a geometric approach,
{\it Ann. Inst. Fourier.}, {\bf 50} (2000), 321-362.

\bibitem{constantin}
A. Constantin, On the Blow-up of solutions of a periodic shallow
water equation, {\it J. Nonlinear Sci.}, {\bf 10} (2000), 391-399.

\bibitem{c}A. Constantin, The trajectories of particles in Stokes waves,
{\it Invent. Math.}, {\bf 166} (2006), 523-535.

\bibitem{ce}A. Constantin and J. Escher, Global existence and blow-up for a shallow water equation,
{\it Annali Sc. Norm. Sup. Pisa.}, {\bf 26} (1998), 303-328.

\bibitem{c-e}A. Constantin and J. Escher, Wave breaking for nonlinear nonlocal shallow water equations,
{\it Acta Mathematica.}, {\bf 181} (1998), 229-243.

\bibitem{ce1}A. Constantin and J. Escher, Global weak solutions for a shallow water equation, {\it Indiana J.
Math.}, {\bf 47} (1998), 1527-1545.

\bibitem{c-j}A. Constantin and R. S. Johnson, Propagation
of very long water waves, with vorticity, over variable depth, with
applications to tsunamis, {\it Fluid Dynam. Res.}, {\bf 40} (2008), 175-211.

\bibitem{c-k-k-t}A. Constantin, T. Kappeler, B. Kolev
and P. Topalov, On geodesic exponential maps of the Virasoro group, {\it Ann. Global Anal. Geom.},
{\bf 31} (2007), 155-180.

\bibitem{C-K}A. Constantin and B. Kolev, On the geometric approach to the motion of inertial mechanical
systems, {\it J. Phys. A}, {\bf 35} (2002), R51-R79.

\bibitem{C-L}A. Constantin and D. Lannes, The hydrodynamical relevance of the Camassa-Holm
and Degasperis-Procesi equations, {\it Arch. Ration. Mech. Anal.}, {\bf 192} (2009), 165-186.

\bibitem{C-M}A. Constantin and H. P. McKean, A shallow water equation on the circle, {\it Comm. Pure
Appl. Math.}, {\bf 52} (1999), 949-982.

\bibitem{cm}A. Constantin and L. Molinet, Obtital stability of solitary waves for a shallow water
equation, {\it Phys. D}, {\bf 157} (2001), 75-89.

\bibitem{cs}A. Constantin and W. A. Strauss, Stability of peakons, {\it Comm. Pure Appl. Math.}, {\bf 53} (2000),
603-610.

\bibitem{HHD}H. H. Dai and M. Pavlov, Transformations for the Camassa-Holm equation, its high-frequency
limit and the Sinh-Gordon equation, {\it J. P. Soc. Japan}, {\bf 67} (1998), 3655-3657.

\bibitem{dm}K. E. Dika and L. Molinet, Stability of multi antipeakon-peakons profile, {\it Disc. Cont. Dyn.
Sys-Ser.B}, {\bf 12} (2009), 561-577.

\bibitem{dm1}K. E. Dika and L. Molinet, Stability of multipeakons, {\it Ann. I. H. Poincar\'{e}}, {\bf 26} (2009),
1517-1532.

\bibitem{escher}
J. Escher, M. Kohlmann and B. Kolev, Geometric aspects of the periodic $\mu$DP equation, 2010. URLhttp://arxiv.org/abs/1004.0978v1.

\bibitem{ewy}J. Escher, S. Wu, and Z. Yin, Global existence and blow-up phenomena for
a weakly dissipative Degasperis-Procesi equation. {\it Discrete Contin. Dyn. Syst. Ser. B}, {\bf 12}, no. 3, (2009), 633-645.

\bibitem{F-F}A. Fokas and B. Fuchssteiner, Symplectic structures, their B$\ddot{a}$cklund transformations and
hereditary symmetries, {\it Phys. D}, {\bf 4} (1981), 47-66.

\bibitem{fu} Y. Fu, Y. Liu and C. Qu, On the blow up structure for the generalized periodic Camassa-Holm and
Degasperis-Procesi equations, arXiv:1009.2466.

\bibitem{g}J. M. Ghidaglia, Weakly damped forced Korteweg-de Vries equations behave
as finite dimensional dynamical system in the long time. {\it J. Differential Equations},
{\bf 74}, (1988), 369-390.

\bibitem{gui}G. L. Gui, Y. Liu and M. Zhu, On the wave-breaking phenomena and global existence for
the generalized periodic Camassa-Holm equation, arXiv:1107.3191v2.

\bibitem{J-R} J. K. Hunter and R. Saxton, Dynamics of director fields, {\it SIAM J. Appl. Math.}, {\bf51} (1991),
1498-1521.

\bibitem{J-Y} J. K. Hunter and Y. Zheng, On a completely integrable nonlinear
hyperbolic variational equation, {\it Phys. D}, {\bf 79} (1994), 361-386.

\bibitem{J-H}R. S. Johnson, Camassa-Holm, Korteweg-de Vries and related models for water waves, {\it J.
Fluid. Mech.}, {\bf 455} (2002), 63-82.

\bibitem{Kato 1}
T. Kato, Quasi-linear equations of evolution, with applications to
partial differential equations, in "Spectral Theory and Differential
Equations", Lecture Notes in Math., Vol. 448, Springer Verlag,
Berlin, (1975), 25--70.

\bibitem{Kato 3}
T. Kato and G. Ponce, Commutator estimates and Navier-Stokes
equations, {\it Comm. Pure Appl. Math.}, {\bf 41} (1988), 203-208.

\bibitem{k-l-m}
B. Khesin, J. Lenells and G. Misiolek, Generalized Hunter-Saxton
equation and the geometry of the group of circle diffeomorphisms,
{\it Math. Ann.}, {\bf 342} (2008), 617-656.

\bibitem{k2}M. Kohlmann, Global existence and blow-up for a weakly diddipative $\mu$DP equation, arXiv:1010.2355v2.

\bibitem{k}B. Kolev, Poisson brackets in hydrodynamics, {\it Discrete Contin. Dyn. Syst.}, {\bf 19} (2007), 555-574.

\bibitem{d}D. Kruse, Variational derivation of the Camassa-Holm shallow
water equation, {\it J. Nonlinear Math. Phys.}, {\bf 14} (2007), 311-320.

\bibitem{l}M. Lakshmanan, Integrable nonlinear wave
equations and possible connections to tsunami dynamics, in
"Tsunami and nonlinear waves", pp. 31-49, Springer, Berlin, 2007.

\bibitem{l1}J. Lenells, The Hunter-Saxton equation describes the geodesic flow on a sphere, {\it J. Geom.
Phys.}, {\bf 57} (2007), 2049-2064.

\bibitem{l-m-t}
J. Lenells, G. Misiolek and F. Ti$\breve{g}$lay, Integrable
evolution equations on spaces of tensor densities and their peakon
solutions, {\it Commun. Math. Phys.}, {\bf 299} (2010), 129-161.

\bibitem{lo}Y. Li and P. Olver, Well-posedness and blow-up solutions for an integrable nonlinearly dispersive
model wave equation, {\it J. Differential Equations}, {\bf 162} (2000), 27-63.

\bibitem{m}G. Misiolek, Classical solutions of the periodic Camassa-Holm equation, {\it Geom. Funct. Anal.},
{\bf 12} (2002), 1080-1104.

\bibitem{P-P} P. Olver and P. Rosenau, Tri-Hamiltonian duality between solitons and solitary wave solutions
having compact support, {\it Phys. Rev. E(3)}, {\bf53} (1996), 1900-1906.

\bibitem{os}E. Ott and R. N. Sudan, Damping of solitary waves. {\it Phys. Fluids}, {\bf 13}, (1970), 1432-1434.

\bibitem{t}J. F. Toland, Stokes waves, {\it Topol. Methods Nonlinear Anal.}, {\bf 7} (1996), 1-48.

\bibitem{w}G. B. Whitham, "Linear and nonlinear waves", Wiley-Interscience, New York-London-Sydney, 1974.

\bibitem{wy}S. Wu and Z. Yin, Blow-up and decay of the solution of the weakly dissipative
Degasperis-Procesi equation. {\it SIAM J. Math. Anal.}, {\bf 40}, no. 2, (2008), 475-490.

\bibitem{wy1}S. Wu and Z. Yin, Global existence and blow-up phenomena for the weakly
dissipative Camassa-Holm equation. {\it J. Differential Equations}, {\bf 246}, no. 11, (2009),  4309-4321.

\bibitem{xz}Z. Xin and P. Zhang, On the weak solutions to a shallow water equation, {\it Comm. Pure Appl.
Math.}, {\bf 53} (2000), 1411-1433.

\bibitem{y}Z. Yin, Well-posedness, global existence and blowup phenomena for an integrable shallow
water equation, {\it Discrete Continuous Dynam. Systems}, {\bf 10} (2004), 393-411.

\bibitem{y1}Z. Yin, Blow-up phenomena and decay for the periodic Degasperis-Procesi
equation with weak dissipation. {\it J. Nonlinear Math. Phys.}, {\bf 15}, (2008),  28-49.

\bibitem{y2}Z. Yin, On the Cauchy problem for the generalized Camassa-Holm equation. {\it Nonlinear Anal.}, {\bf 66}, (2007), 460-471.

\bibitem{Y} Z. Yin, On the structure of solutions to the periodic Hunter-Saxton equation, {\it SIAM J. Math.
Anal.}, {\bf 36} (2004), 272-283.

\bibitem{Y1} Z. Yin, Global existence for a new periodic integrable
equation, {\it J. Math. Anal. Appl.}, {\bf 283}, (2003), 129-139.

\end{thebibliography}
\end{document}